\documentclass{amsart}

\usepackage[xdvi]{graphics}
\usepackage{amsthm}
\usepackage{amssymb}

\newtheorem{theorem}{Theorem}[section]
\newtheorem{corollary}[theorem]{Corollary}
\newtheorem*{proposition*}{Proposition}

\newtheorem*{theorem*}{Theorem}
\newtheorem{lemma}[theorem]{Lemma}

\newtheorem{claim}[theorem]{Claim}
\newtheorem{condition}[theorem]{Condition}

\theoremstyle{definition}

\theoremstyle{remark}
\newtheorem{remark}[theorem]{Remark}

\newcommand{\Diagram}{\mathcal{D}}
\newcommand{\BasicEdgepath}{\lambda}
\newcommand{\BasicEdgepathSystem}{\Lambda}
\newcommand{\NearlyBasicEdgepath}{\lambda^\prime}

\newcommand{\Edgepath}{\gamma}
\newcommand{\EdgepathSystem}{\Gamma}
\newcommand{\angleb}[1]{\langle #1 \rangle} 
\newcommand{\circleb}[1]{\langle #1 \rangle^{\circ}} 

\newcommand{\oriented}[1]{\protect \overrightarrow{#1}}
\newcommand{\notoriented}[1]{#1}
\newcommand{\Vinf}{V_\infty}
\newcommand{\Vone}{V_1}
\newcommand{\Dinf}{D_\infty}
\newcommand{\Dzero}{D_0}
\newcommand{\Hone}{H_1}
\newcommand{\Hzero}{H_0}

\newcommand{\Diam}{\mathrm{Diam}}

\newcommand{\VDiagram}{
           \put(0,0){\scalebox{0.3}{\includegraphics{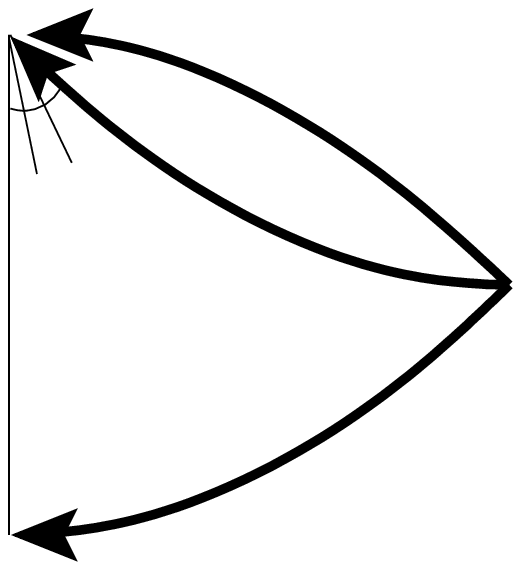}}}
           \put(-42,0){
            \put(60,60){$\NearlyBasicEdgepath_s(\oriented{T_A})$}
            \put(70,56){\vector(0,-1){25}}
            \put(36,52){$\angleb{1}$}
            \put(36,-7){$\angleb{0}$}
            \put(75,38){$\BasicEdgepath_\inc(\notoriented{T_A})$}
            \put(68,3){$\BasicEdgepath_\dec(\notoriented{T_A})$}
            \put(53,25){\vector(-1,2){5}}
            \put(45,17){odd}
           }
}

\newcommand{\Twist}{{\rm twist}}
\newcommand{\Slope}{{\rm slope}} 
\newcommand{\slope}{R}

\newcommand{\Crossing}{{\rm C}}

\newcommand{\dec}{{\rm dec}}
\newcommand{\inc}{{\rm inc}}
\newcommand{\rmmin}{{\rm min}}
\newcommand{\rmmax}{{\rm max}}

\newcommand{\NumTangles}{N}
\newcommand{\TangleP}{P}
\newcommand{\TangleQ}{Q}
\newcommand{\TangleR}{R}
\newcommand{\TangleS}{S}

\begin{document}
\title
[%
Boundary slopes and the numbers of positive/negative crossings
]
{%
Boundary slopes and the numbers of positive/negative crossings for Montesinos knots
}

\author{Kazuhiro Ichihara}
\address{%
School of Mathematics Education,
Nara University of Education,
Takabatake-cho, Nara 630--8528, Japan}
\email{ichihara@nara-edu.ac.jp}

\author{%
    Shigeru Mizushima}
\address{%
        Department of Mathematical and Computing Sciences \\
        Tokyo Institute of Technology \\
        2--12--1 Ohokayama, Meguro \\
        Tokyo 152--8552, Japan}
\email{mizusima@is.titech.ac.jp}

\thanks{The first author is partially supported by 
Grant-in-Aid for Young Scientists (B), No. 20740039, 
Ministry of Education,Culture,Sports,Science and Technology, Japan.}


\begin{abstract}
We show that
a finite numerical boundary slope of an essential surface in the exterior of a Montesinos knot
is bounded above and below in terms of the numbers of positive/negative crossings
of a specific minimal diagram of the knot.
\end{abstract}

\maketitle


\section{Introduction}

Let $K$ be a knot in the $3$-sphere $S^3$ with the exterior $E(K)$.
The boundary of a properly embedded 
essential (i.e., incompressible and boundary-incompressible) 
surface in $E(K)$ gives 
a parallel family of non-trivial simple closed curves 
on the boundary torus $\partial E(K)$. 
It is expressed by an irreducible fraction (possibly infinity $1/0$), 
which is called a {\em numerical boundary slope}. 
See \cite{R} for example. 
%
%
We say that a boundary slope other than $1/0$ is {\em finite}.

%
In this paper, 
we consider the knot $K$ called a \textit{Montesinos knot} in $S^3$, 
which is a knot obtained by connecting rational tangles. 
See the next section for detail. 
For Montesinos knots, 
there exists an algorithm developed by Hatcher and Oertel 
which completely determines the set of boundary slopes (\cite{HO}).
Based on their algorithm, 
we present an estimate on 
the concrete range of boundary slope as follows.
%
\begin{theorem}
\label{MainThm}
\label{Thm:Main}
Let $K$ be a Montesinos knot 
consisting of rational tangles $\TangleP_1/\TangleQ_1$, $\TangleP_2/\TangleQ_2$, $\ldots$, $\TangleP_\NumTangles/\TangleQ_\NumTangles$.
Assume that the number $\NumTangles$ of the rational tangles is three or more, and each fraction $\TangleP_i/\TangleQ_i$ is a non-integral non-infinity fraction.
Then,
we have
\[
-2\,\Crossing_-(D)\le \slope \le 2\,\Crossing_+(D)
\]
for any finite boundary slope $\slope$ for $K$ and 
any standard diagram $D$ of $K$, 
where $\Crossing_+(D)$ and $\Crossing_-(D)$
denote the number of positive crossings and the number of negative crossings
of the diagram $D$ respectively.
\end{theorem}
%

\begin{remark}\label{Rem:StandardMinimal}
A \textit{standard diagram} $D$ is a specific diagram 
naturally obtained for a Montesinos knot.
See the next section for detail. 
Actually, we can choose a standard diagram $D$ so that $D$ attains the crossing number.
\end{remark}

\subsection{Corollary}

It is known that 
at least two distinct numerical boundary slopes always exist for a non-trivial knot \cite{CS84},
and the number of such slopes is finite for any knot \cite{H}.
Hence, the set of boundary slopes 
gives a nonempty finite subset of $\mathbb{Q}\cup \{1/0\}$.
This is called the {\em (numerical) boundary slope set}.
%
%
We define the {\em diameter} of the boundary slope set
as the maximum minus the minimum among finite boundary slopes.
Note that the infinity slope is ignored in this definition.

%
In \cite{IMCrossing}, 
we obtained an estimate of the diameter by the crossing number
for a Montesinos knot $K$ as follows:
\begin{eqnarray}
\Diam(K)\le 2\,\Crossing(K)
 \label{Eq:IneqDiameter}
\end{eqnarray}
where $\Diam(K)$ is the diameter of the boundary slope set for $K$
and $\Crossing(K)$ is the minimal crossing number.
%
%
Now this inequality (\ref{Eq:IneqDiameter})
is easily obtained from Theorem \ref{MainThm}.
Moreover,
together with \cite{MMR},
it becomes a corollary of Theorem \ref{MainThm}
as follows: 
\begin{corollary}
For a Montesinos knot $K$, we have the inequality
\[
\Diam(K)\le 2\,\Crossing(K).
\]
\end{corollary}
\begin{proof}
For 2-bridge knots, the result is shown in \cite{MMR}.
Assume now that the Montesinos knot $K$ is not a 2-bridge knot. 
We consider the expression of $K$ 
so that its standard diagram attains the minimal crossing number. 
Then, we have
$\Diam(K)\le 2\,\Crossing_{+}(D)-(-2\,\Crossing_{-}(D))=2\,(\Crossing_{+}(D)+\Crossing_{-}(D))=2\,\Crossing(D)=2\,\Crossing(K)$.
\end{proof}
%

\subsection{Motivations and related results}

Our study is originally motivated by the unpublished work \cite{IS} 
on the boundary slopes and the crossing number of knots. 
Suggested by their work, we obtained the corollary above in \cite{IMCrossing}. 

Extending our work in \cite{IMCrossing}, it was shown in \cite{MMR} that 
Inequality (\ref{Eq:IneqDiameter}) holds for 2-bridge knots, 
and also shown that a similar inequality holds for 2-bridge link in \cite{HoSha}. 

Also motivated by \cite{IS} and \cite{IMCrossing}, 
Yamaoka \cite{Y} obtained 
\[
 \slope_\rmmax= 2\,\Crossing_+(D_a), 
 \text{\ \ and \ \ } 
 \slope_\rmmin= -2\,\Crossing_-(D_a), 
\]
for a reduced alternating diagram $D_a$ of a two-bridge knot $K$, 
where $\slope_\rmmax$ and $\slope_\rmmin$ denote
the maximum and the minimum among the finite boundary slopes.
Thus, for a two-bridge knot $K$, we have
\[
-2\,\Crossing_-(D_a)\le \slope \le 2\,\Crossing_+(D_a) \ \ .
\]

Actually, this equality holds for arbitrary diagram for a two-bridge knot. 
Because it was observed in \cite{IS} that 
$\Crossing_+(D_a) \le \Crossing_+(D)$ and $\Crossing_-(D) \le \Crossing_-(D_a)$ 
for a reduced alternating diagram $D_a$ and any diagram $D$ of an alternating knot. 
In fact, this follows from \cite[Theorem 13.3]{M} immediately. 

In the case of general alternating knots, 
Hayashi observed that, 
for a reduced alternating diagram $D_a$, 
$2\,\Crossing_+(D_a)$ and $-2\,\Crossing_-(D_a)$ are 
equal to the boundary slopes of checkerboard surfaces \cite{Ha}. 
Thus if the boundary slopes of checkerboard surfaces 
attain the maximal and minimal boundary slopes, 
then we have the above equalities immediately. 

For example, for alternating Montesinos knots, 
we showed in \cite{IMCrossing} 
that the boundary slopes of checkerboard surfaces 
attain the maximal and minimal boundary slopes. 
And so, for alternating Montesinos knots, 
we have the following for any diagram $D$; 
\[
-2\,\Crossing_-(D)\le \slope \le 2\,\Crossing_+(D)
\]

%

%
%
Also note that, in Theorem \ref{Thm:Main},
the inequality holds just for specific diagrams 
for Montesinos knots including non-alternating ones.
It might be a supporting evidence
of the inequality
for non-alternating Montesinos knots and its arbitrary diagram.

About Inequality (\ref{Eq:IneqDiameter}), 
it is natural to ask whether it holds for other class of knots. 
Actually it holds for all torus knots suggested in \cite{MMR}, 
and for all known example until recently. 
However, very recently, Kabaya gave \cite{K} examples of knots 
which have boundary slope not satisfying the inequality in a completely different method.


\bigskip

This paper is organized as follows.
We review some definitions and
Hatcher-Oertel algorithm in Section 2,
and then prove the main theorem in Section 3.



\section{Preliminary}

%
The readers familiar to \cite{HO}, 
or \cite{IMBounds, IMCrossing,IMLowerBound} 
would be able to skip this section. 

\begin{figure}[hbt]
 \begin{picture}(140,65)
  \put(0,8){\scalebox{0.6}{\includegraphics{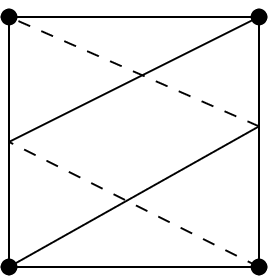}}}
  \put(80,0){\scalebox{0.6}{\includegraphics{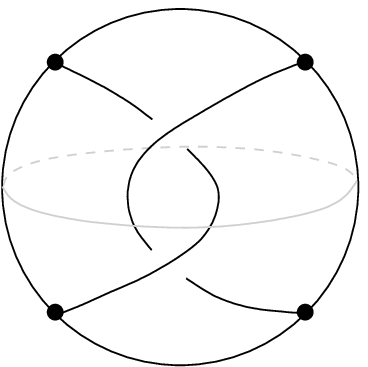}}}
 \end{picture}
 \caption{Arcs on a fourth-punctured sphere and $1/2$-tangle}
 \label{Fig:1over2-Tangle}
\end{figure}

\subsection{Montesinos knot}

\ \\
%
%
\noindent {\bf Rational tangle}\\
Let $\TangleP/\TangleQ$ be an irreducible fraction.
We draw segments of slope $\TangleP/\TangleQ$ on the flattened fourth-punctured sphere. The segments form two arcs on the sphere.
Then we push the arcs into the interior of the sphere
by an isotopy fixing endpoints of arcs.
The two-string tangle thus obtained is called a {\em rational tangle}.
In particular, a rational tangle corresponding to an irreducible fraction $\TangleP/\TangleQ$ is called a {\em $\TangleP/\TangleQ$-tangle}.
See Figure \ref{Fig:1over2-Tangle}.
In this paper,
we may use the term ``tangle'' as ``rational tangle'' for ease.

%
%

\begin{figure}[htb]
 \begin{center}
  \begin{picture}(109,91)
   \put(0,0){\scalebox{0.25}{\includegraphics{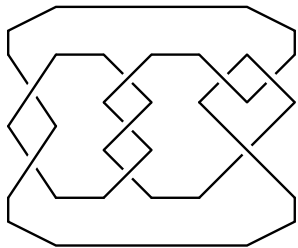}}}
  \end{picture}
  \caption{A diagram of $M(1/2,1/3,-2/3)$}
  \label{Fig:MontesinosKnot}
 \end{center}
\end{figure}

%
%
\noindent {\bf Montesinos knot}\\
For irreducible fractions $\TangleP_i/\TangleQ_i$ $(i=1,2,\ldots,N)$,
a {\em Montesinos knot} $M(\TangleP_1/\TangleQ_1$, $\TangleP_2/\TangleQ_2$, $\ldots$, $\TangleP_\NumTangles/\TangleQ_\NumTangles)$ is a knot
obtained by
connecting rational tangles corresponding to $\TangleP_i/\TangleQ_i$'s.
See Figure \ref{Fig:MontesinosKnot}.
Throughout this paper,
we fix some notations and assumptions for a Montesinos knot.
$\NumTangles$ denotes the number of tangles of a Montesinos knot.
Moreover, $\TangleP_i/\TangleQ_i$ denotes a fraction 
representing the $i$-th tangle $T_i$ in the Montesinos knot.
According to \cite{HO},
we assume that each $\TangleP_i/\TangleQ_i$ is not an integer, and is not $1/0$.
Assume further that $N\ge 3$.
With these conditions,
a Montesinos knot is normalized in a sense,
and then two-bridge knots and non-prime knots are excluded from argument.

For a Montesinos link to be a knot,
or equivalently,
to have only one link component,
the tuple $(\TangleP_1/\TangleQ_1, \TangleP_2/\TangleQ_2, \ldots, \TangleP_\NumTangles/\TangleQ_\NumTangles)$ of fractions
must satisfy either of the following two conditions.
\begin{condition}
 \label{Condition:MontesinosKnot}\mbox{}\\ \vspace{-6mm}
\begin{itemize}
 \item[(a)]
       Exactly one fraction has even denominator.
       Any other fraction has odd denominator.
 \item[(b)]
       All fraction has odd denominator.
       Furthermore,
       the number of fractions with odd numerator
       must be odd.
\end{itemize}
\end{condition}

Usually, two isotopic knots are identified with each other.
In contrast,
an expression $M(\TangleP_1/\TangleQ_1,\TangleP_2/\TangleQ_2,\ldots,\TangleP_\NumTangles/\TangleQ_\NumTangles)$ of a Montesinos knot
may be thought to specify 
not only an isotopy class of a knot
but also a decomposition of a knot into a tuple of rational tangles.
We adopt this point of view in this paper.
Note that one isotopy class of Montesinos knots is
actually expressed by many different expressions as a Montesinos knot.

\medskip

%
%
\noindent {\bf Standard diagram of a rational tangle}\\
A fixed fraction $\TangleP/\TangleQ$ can be expressed as a {\em standard continued fraction} 
\[
\frac{\TangleP}{\TangleQ}=
a_1+
\frac{1}{a_2+
\frac{1}{\ddots
\frac{1}{a_{k-1}+\frac{1}{a_k}
}}}
,
\]
where the sequence $[a_1,a_2,\ldots,a_k]$ of integers
satisfies
$k\ge 1$,
$a_1\ge 0$,
$a_i\ge 1$ for $i=2,3,\ldots,k-1$,
and $a_k\ge 2$ if $\TangleP/\TangleQ$ is positive,
and
$k\ge 1$,
$a_1\le 0$,
$a_i\le -1$ for $i=2,3,\ldots,k-1$,
and $a_k\le -2$ if $\TangleP/\TangleQ$ is negative.
We call $k$ the {\em length} of the expansion.
%
%
%
%
%

\begin{figure}[hbt]
        \begin{center}
         \begin{picture}(380,80)
          \put(  0,0){
           \put(0,0){\scalebox{0.8}{\includegraphics{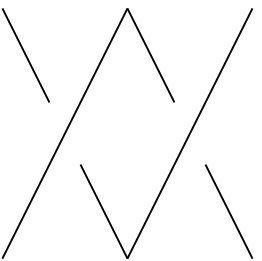}}}
           \put(10,-15){$2$-tangle}
           \put(90,0){\scalebox{0.8}{\includegraphics{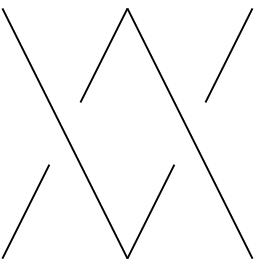}}}
           \put(100,-15){$-2$-tangle}
          }
          \put(190,0){
           \put(0,0){\scalebox{0.3}{\includegraphics{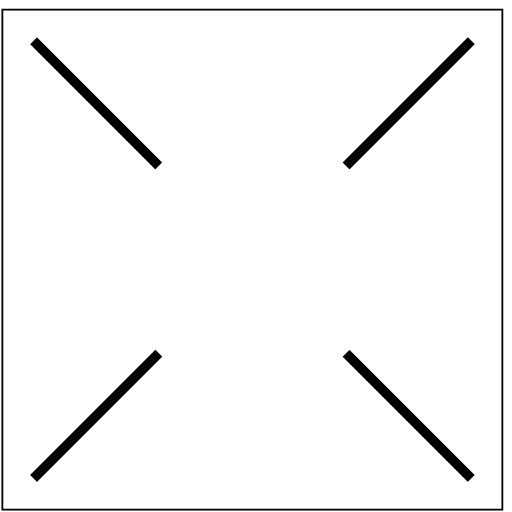}}}
           \put(16,18){$D_A$}
          }
          \put(265,0){
           \put(0,0){\scalebox{0.3}{\includegraphics{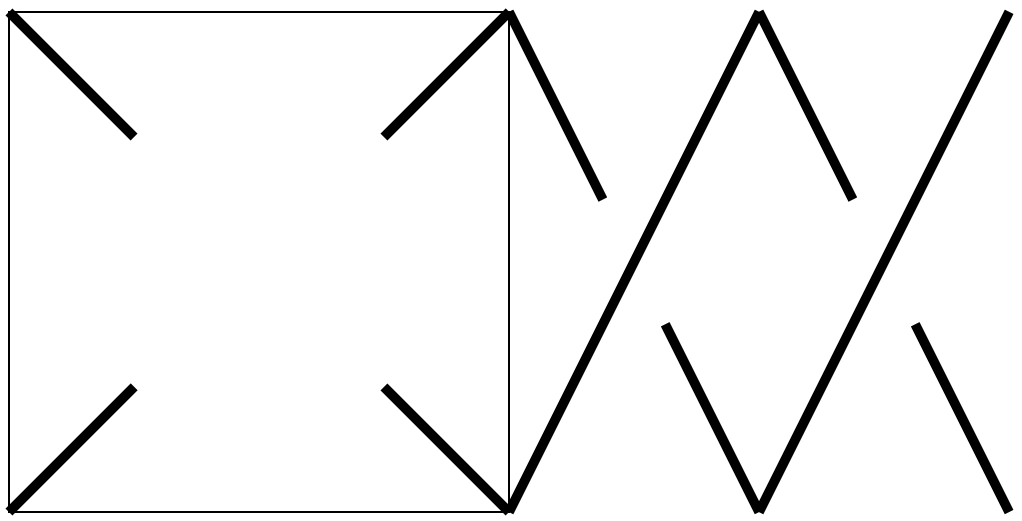}}}
           \put(18,16){\rotatebox{90}{\reflectbox{$D_A$}}}
           \put(48,50){\scalebox{1.0}[1.0]{$\overbrace{\hspace{14mm}}$}}
           \put(50,61){$a$-tangle}
          }
         \end{picture}
        \end{center}
        \caption{standard diagrams}
        \label{Fig:StandardDiagrams}
       \end{figure}

According to a fraction expansion,
we can make a diagram of $\TangleP/\TangleQ$-tangle.
For a tangle corresponding to $[a]=a$,
the diagram is like the left two diagrams of Figure \ref{Fig:StandardDiagrams}.
That is, the diagram is made by aligning crossings horizontally.
Next, 
if $D_A$ is the diagram of a rational tangle $T_A$
corresponding to a fraction expansion $[a_1,a_2,\ldots,a_k]$,
then 
the diagram of some rational tangle $T_B$
corresponding to a fraction expansion $[a,a_1,a_2,\ldots,a_k]$
is obtained by combining 
a mirror image of $D_A$ and 
$a$-tangle
where the mirror image is taken with respect to a line from top left to bottom right
(see the right diagram of Figure \ref{Fig:StandardDiagrams}). 
Inductively, a diagram is defined for a fraction expansion of $\TangleP/\TangleQ$.
In particular,
the {\em standard diagram} of a rational tangle
is a diagram for the standard fraction expansion of the fraction.
%
%


%
%
\noindent{\bf Standard diagram of a Montesinos knot}\\
A {\em standard diagram of a Montesinos knot}
$K=M(\TangleP_1/\TangleQ_1$, $\TangleP_2/\TangleQ_2$, $\ldots$, $\TangleP_\NumTangles/\TangleQ_\NumTangles)$
is obtained by
collecting standard diagrams of $\TangleP_i/\TangleQ_i$-tangles
and aligning the standard diagrams horizontally in line.
For instance,
Figure \ref{Fig:MontesinosKnot} is the standard diagram 
of the Montesinos knot $M(1/2,1/3,-2/3)$.
Since the standard diagram
is constructed for a tuple of rational tangles,
strictly speaking,
the standard diagram is defined 
not for an isotopy class of a Montesinos knot
but for an expression of a Montesinos knot.

\begin{figure}[hbt]
        \begin{center}
         \begin{picture}(150,60)
          \put(  0,0){
           \put(0,10){\scalebox{1.0}{\includegraphics{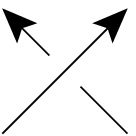}}}
           \put(2,-3){positive}
          }
          \put(100,0){
           \put(0,10){\scalebox{1.0}{\includegraphics{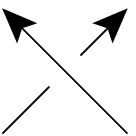}}}
           \put(1,-3){negative}
          }
         \end{picture}
        \end{center}
        \caption{Sign of crossing}
        \label{Fig:SignOfCrossing}
       \end{figure}

%
%
\noindent{\bf The number of positive/negative crossings}\\
%
%
%
%
For a diagram $D$ of a knot $K$,
we give either of two possible orientations.
Then, a sign is given to each crossing,
where
we adopt the convention of signs of crossings
as in Figure \ref{Fig:SignOfCrossing}.
Then,
we can determine the number of positive crossings
and the number of negative crossings,
which are denoted by 
$\Crossing_{+}(D)$ and $\Crossing_{-}(D)$ respectively.
%


\subsection{Hatcher-Oertel algorithm}

Here, we briefly review Hatcher-Oertel algorithm.
We will try to give sufficient explanations 
so that this paper can be self-contained, 
and so, we omit certain part of their machinery. 
Please see \cite{HO}
or our previous papers \cite{IMBounds} or \cite{IMCrossing} for more detail.

The Hatcher-Oertel algorithm determines the set of all numerical boundary slopes of a Montesinos knot.

In the algorithm,
$S^3$ is decomposed into $N$ 3-balls
where
a Montesinos knot is also decomposed into $N$ rational tangles at the same time.
Furthermore,
an essential surface is also decomposed into $N$ pieces of surfaces,
which are called ``subsurfaces''.
A subsurface around a rational tangle
is expressed by an ``edgepath'',
and an essential surface for a Montesinos knot
is expressed by a tuple of $N$ edgepaths,
which is called an ``edgepath system''.

\begin{figure}[htb]
  \begin{center}
   \begin{minipage}{180pt}
   \begin{center}
    \scalebox{0.75}{
    \begin{picture}(70,150)
     \put(0,0){\scalebox{0.7}{\includegraphics{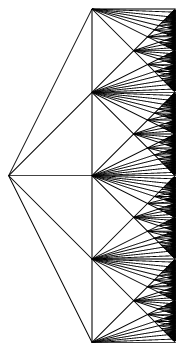}}}
     \put(49,135){\rotatebox{90}{\scalebox{1.0}{$\cdots$}}}
     \put(49,-10){\rotatebox{90}{\scalebox{1.0}{$\cdots$}}}
     \put(2,70.5){\scalebox{1.5}{\vector(1,0){55}}}
     \put(42,5){\scalebox{1.5}{\vector(0,1){90}}}
     \put(29,15){\tiny $-2$}
     \put(29,40){\tiny $-1$}
     \put(36,98){\tiny $1$}
     \put(36,123){\tiny $2$}
     \put(40,148){$v$}
     \put(3,82){\vector(1,-1){10}}
     \put(-10,85){$\angleb{1/0}$}
     \put(2,64){\tiny $-1$}
     \put(72,64){\tiny $1$}
     \put(92,68){$u$}
     \put(33,63){\tiny $O$}
    \end{picture}
    }
   \end{center}
   \caption{The diagram $\Diagram$}
   \label{Fig:Diagram}
   \end{minipage}
\hspace{-10mm}
    \begin{minipage}{200pt}
    \begin{center}    
    \scalebox{0.75}{
     \begin{picture}(150,150)
      \put(0,-10){\scalebox{0.75}{\includegraphics{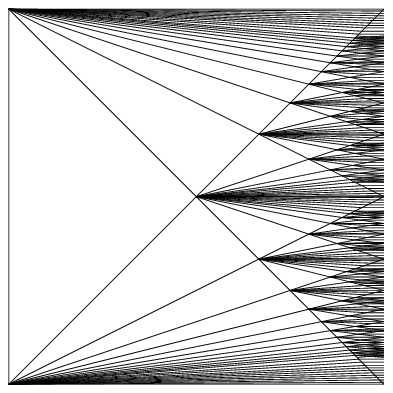}}}
      \put(-5,-2){\vector(1,0){10}}
      \put(-21,-5){$\angleb{0}$}
      \put(156,-2){\vector(-1,0){10}}
      \put(158,-5){$\circleb{0}$}
      \put(-5,133){\vector(1,0){10}}
      \put(-21,130){$\angleb{1}$}
      \put(156,133){\vector(-1,0){10}}
      \put(158,130){$\circleb{1}$}
      \put(63,66){\vector(1,0){10}}
      \put(46,63){$\angleb{\frac{1}{2}}$}
      \put(156,66){\vector(-1,0){10}}
      \put(158,63){$\circleb{\frac{1}{2}}$}
      \put(83,43){\vector(1,0){10}}
      \put(66,40){$\angleb{\frac{1}{3}}$}
      \put(156,43){\vector(-1,0){10}}
      \put(158,40){$\circleb{\frac{1}{3}}$}
      \put(83,88){\vector(1,0){10}}
      \put(66,85){$\angleb{\frac{2}{3}}$}
      \put(156,88){\vector(-1,0){10}}
      \put(158,85){$\circleb{\frac{2}{3}}$}
     \end{picture}
    }
    \end{center}
    \caption{
     A part of the diagram $\Diagram$ in 
     $[0,1]\times[0,1]$
    }
    \label{Fig:Diagram2}
    \end{minipage}

  \end{center}
  \end{figure}

An edgepath is a path in a {\it diagram} $\Diagram$. 
The diagram $\Diagram$ is the 1-skeleton of a triangulation of 
a region $\{(u,v)|-1<u<1\}$ in a $uv$-plane.
See Figure \ref{Fig:Diagram} and Figure \ref{Fig:Diagram2}.
Vertices are points $\angleb{p/q}$ for $p/q \in \mathbb{Q}\cup \{1/0\}$
and points $\circleb{p/q}$ for $p/q \in \mathbb{Q}$.
Two points $\angleb{p/q}$ and $\angleb{r/s}$ are connected by an edge
if $|ps-qr|=1$.
An edge $\angleb{1/0}$\,--\,$\angleb{z}$ for an integer $z$ is called an $\infty$-edge.
An edge $\angleb{z}$\,--\,$\angleb{z+1}$ for an integer $z$ is called a vertical edge.
$\angleb{p/q}$ and $\circleb{p/q}$ are also connected by an edge,
which is called a {\em horizontal edge}.

An edgepath is a path on $\Diagram$
which is a combination of edges of $\Diagram$.
Roughly right-to-left orientation is given to an edgepath.
For a fraction $\TangleP/\TangleQ$,
the starting point of an edgepath for $\TangleP/\TangleQ$-tangle
must lie on the horizontal edge 
$\angleb{\TangleP/\TangleQ}$\,--\,$\circleb{\TangleP/\TangleQ}$.
See \cite{IMCrossing} for concrete examples of edgepaths.
An edgepath is expressed by a sequence of vertices of $\Diagram$.

See \cite{IMBounds} or \cite{IMCrossing}
for the correspondence between an edgepath and a subsurface.
Roughly, an edge corresponds to a saddle
and an edgepath (or an edgepath system) corresponds to a combination of saddles.

%
%

Each edgepath must satisfy ``minimality condition'',
which is required 
so that corresponding subsurface is essential.
An edgepath system must satisfy ``gluing consistency'',
which assures that subsurfaces corresponding to edgepaths in an edgepath system can be consistently glued each other.
By the gluing consistency,
endpoints of all edgepaths in an edgepath system 
have common $u$-coordinate (first coordinate).
According to the common coordinate,
edgepath systems are classified into type I/II/III.
Edgepaths of an edgepath system of type I have endpoints with $u>0$.
Edgepaths of an edgepath system of type II have endpoints with $u=0$,
in short, endpoints lie on the vertical line $u=0$.
Edgepaths of an edgepath system of type III have endpoints with $u<0$.
That is, the endpoints are $\angleb{1/0}$.

A {\em basic edgepath} is an edgepath which ends at an integer vertex $\angleb{z}$ for some integer $z$ 
and does not include vertical edges.
A {\em basic edgepath system} is an edgepath system
which consists of basic edgepaths.
An edgepath system which does not satisfy gluing consistency is said to be {\em formal}.
Generically, a basic edgepath system is formal.

For a non-$\infty$-edge of $\Diagram$,
with a right-to-left orientation,
if $v$-coordinate increases as a point moves from the starting point to the ending point along the edge,
then the edge is said to be {\em increasing}.
Similarly, if $v$-coordinate decreases, the edge is said to be {\em decreasing}.

In Hatcher-Oertel algorithm,
for a Montesinos knot
$M(\TangleP_1/\TangleQ_1$, $\TangleP_2/\TangleQ_2$, $\ldots$, $\TangleP_\NumTangles/\TangleQ_\NumTangles)$,
all possible edgepath systems satisfying 
minimality and gluing consistency
are collected.
%
Then,
each edgepath system is checked 
whether its corresponding surface is essential or not
by the detailed conditions described in \cite{HO}.
Eventually,
the numerical boundary slopes for an edgepath system
corresponding to an essential surface
is calculated.
By collecting the numerical boundary slopes,
we obtain the set of numerical boundary slopes.

The boundary slope is calculated from an edgepath system roughly as follows.
We first calculate ``twist'' of an edgepath system.
The twist of an edge $e$ is defined as $+2$ or $-2$
according to whether the edge is a decreasing or increasing leftward edge respectively.
The twist 
of an edgepath $\Edgepath$ or an edgepath system $\EdgepathSystem$ is the sum of twists of edges in the edgepath or edgepath system.
Let $\Slope(\EdgepathSystem)$ and $\Twist(\EdgepathSystem)$ denote 
the boundary slope and the twist of $\EdgepathSystem$.
The boundary slope of an edgepath system $\EdgepathSystem$
is calculated by
$\Slope(\EdgepathSystem)=\Twist(\EdgepathSystem)-\Twist(\EdgepathSystem_s)$
where $\EdgepathSystem_s$ is the edgepath system of a Seifert surface for a Montesinos knot.




\section{Proof of the main theorem}

In this section, we will give a proof of our main theorem. 
We will prove three key lemmas, 
Lemma \ref{Lem:TwistInequality}, \ref{Lem:SeifertEdgepathSystem}, and \ref{Lem:IdentitiesForKnot}, 
after preparing some terminologies, 
in subsections \ref{subsection31}, \ref{subsection32}, and \ref{subsection33}, respectively. 
Combining these, in subsections \ref{subsection34}, we will give a proof of our main theorem. 
In subsections \ref{subsection35}, we will consider about Remark \ref{Rem:StandardMinimal}.


%
%



\subsection{Monotonic edgepath systems and bounds of boundary slopes}\label{subsection31}


In this subsection,
we define
a monotonically increasing edgepath system $\EdgepathSystem_\inc$
and
a monotonically decreasing edgepath system $\EdgepathSystem_\dec$
for a Montesinos knot $K$.
Then,
these edgepath systems give bounds of boundary slopes.
These bounds have been essentially given in \cite{IMCrossing} or \cite{IMLowerBound}.

\noindent {\bf Monotonic basic edgepath systems}\\
%
%
A basic edgepath is called {\em monotonically increasing}
if all edges in the edgepath are increasing.
A basic edgepath system is called
{\em monotonically increasing}
if each edgepath in the edgepath system 
is monotonically increasing.
The term {\em monotonically decreasing}
is defined similarly
for a basic edgepath and a basic edgepath system.
%
Both the monotonically increasing basic edgepath system
and the monotonically decreasing basic edgepath system
are unique for a fixed Montesinos knot.
Let
$\BasicEdgepathSystem_\inc=(\BasicEdgepath_{\inc,i})$
and 
$\BasicEdgepathSystem_\dec=(\BasicEdgepath_{\dec,i})$
denote the monotonic basic edgepath systems respectively.

\noindent {\bf Monotonic edgepath systems}\\
Next,
from monotonic basic edgepaths,
we define
a monotonically increasing edgepath $\Edgepath_\inc$
and
a monotonically decreasing edgepath $\Edgepath_\dec$
for a $\TangleP/\TangleQ$-tangle as follows.
If $\TangleP/\TangleQ$ is positive,
$\Edgepath_\inc$ is a monotonically increasing type III edgepath
and 
$\Edgepath_\dec$ is a monotonically decreasing type II edgepath
obtained from a monotonically decreasing basic edgepath 
$\BasicEdgepath_\dec$
by adding downward vertical edges connecting
$\angleb{v_{0} (\BasicEdgepath_\dec)}$
and $\angleb{0}$,
where
$v_{0} (\BasicEdgepath_\dec)$ denotes
the $v$-coordinate of the intersection point
between
the basic edgepath $\BasicEdgepath_\dec$ and the vertical axis $u=0$,
or equivalently
the integer representing the left ending vertex of $\BasicEdgepath_\dec$.
If $\TangleP/\TangleQ$ is negative,
$\Edgepath_\dec$ is a monotonically decreasing type III edgepath
and
$\Edgepath_\inc$ is a monotonically increasing type II edgepath
obtained from a monotonically increasing basic edgepath
$\BasicEdgepath_\inc$
by adding upward vertical edges connecting
$\angleb{v_{0} (\BasicEdgepath_\inc)}$
and $\angleb{0}$.


For a Montesinos knot $K=M(\TangleP_1/\TangleQ_1,\TangleP_2/\TangleQ_2,\ldots,\TangleP_\NumTangles/\TangleQ_\NumTangles)$,
$\EdgepathSystem_\inc=(\Edgepath_{\inc,i})$
and 
$\EdgepathSystem_\dec=(\Edgepath_{\dec,i})$
are obtained by
collecting $\Edgepath_\inc$ and $\Edgepath_\dec$
defined as above for each $\TangleP_i/\TangleQ_i$-tangle.
$\EdgepathSystem_\inc$ and $\EdgepathSystem_\dec$ may contain both type II and type III edgepaths and hence may be formal.
However, 
it does not matter,
since we only use twists of these edgepath systems.

%
%
%

\noindent {\bf A bound of boundary slopes} \\
The twists or boundary slopes of $\EdgepathSystem_\inc$ and $\EdgepathSystem_\dec$
%
%
give lower and upper bounds of twists or boundary slopes of essential surfaces.
%
%
That is:
\begin{lemma}
\label{Lem:TwistInequality}
Let $K$ be a Montesinos knot.
Let $\EdgepathSystem$ 
be an edgepath system
corresponding to an essential surface for the knot $K$
with finite boundary slope.
Then,
we have
\begin{eqnarray}
\label{Eq:UpperLowerBound:Twist}
&& \Twist(\EdgepathSystem_\inc) \le
\Twist(\EdgepathSystem) \le \Twist(\EdgepathSystem_\dec)
,
\end{eqnarray}
or equivalently,
\begin{eqnarray}
\label{Eq:UpperLowerBound}
&& \Slope(\EdgepathSystem_\inc) \le
\Slope(\EdgepathSystem) \le \Slope(\EdgepathSystem_\dec)
.
\end{eqnarray}
\end{lemma}
\begin{proof}
If $\EdgepathSystem$ is a type I edgepath system or a type III edgepath system,
by Proposition 3.2 in \cite{IMCrossing},
we have 
$\Twist(\BasicEdgepathSystem_\inc)\le\Twist(\EdgepathSystem)\le\Twist(\BasicEdgepathSystem_\dec)$.
This inequality gives the inequality (\ref{Eq:UpperLowerBound:Twist}). 
If $\EdgepathSystem$ is a type II edgepath system,
with a care for vertical edges in $\EdgepathSystem$,
we can also show the inequality (\ref{Eq:UpperLowerBound:Twist}) for $\EdgepathSystem$.

Since the boundary slope of an edgepath system
is calculated by
the twist of the edgepath system
minus
the twist of a Seifert surface,
the inequality (\ref{Eq:UpperLowerBound}) is easily obtained from
the inequality (\ref{Eq:UpperLowerBound:Twist}).
\end{proof}


\subsection{Seifert surface and its edgepath system}\label{subsection32}

In order to calculate the boundary slope of an edgepath system,
we have to know also the twist of the edgepath system of a Seifert surface.
Therefore,
we investigate the edgepath system of a Seifert surface.
In this part,
we define an edgepath system $\EdgepathSystem_s$ 
and confirm that $\EdgepathSystem_s$ is the edgepath system of a Seifert surface.

\noindent {\bf Oriented tangle} \\
As a preparation, we first define oriented rational tangles.
For a rational tangle,
we can give an orientation for each of two strings of the tangle.
We call it an {\em oriented tangle}.
Next,
we identify
one oriented tangle $\oriented{T_A}$
with 
another oriented tangle $\oriented{T_B}$
obtained from the tangle $\oriented{T_A}$
by reversing the orientations of the two strings at the same time.
Throughout this paper, oriented tangles are always handled with this identification.
Of course, results in this paper are consistent under this identification.
Note that,
under the identification,
one tangle has two choices of orientation.

\begin{figure}[hbt]
 \begin{center}
  \begin{picture}(310,110)
   \put( 50,80){
    \put(-43,13){$\Hzero=$}
    \put(-15,8){\scalebox{1.5}[2.8]{$\left\{\right.$}}
    \put( 0,0){\scalebox{0.2}{\includegraphics{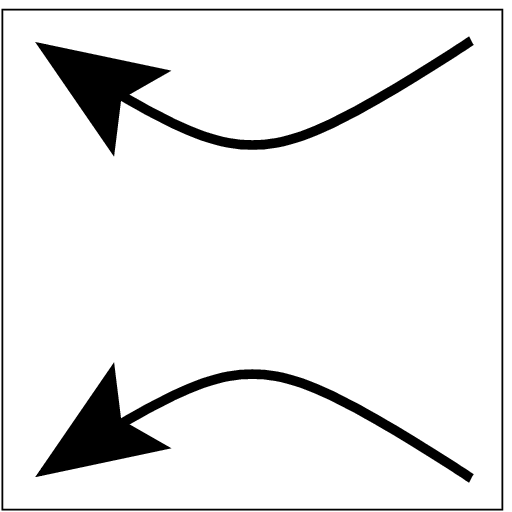}}}
    \put(38,3){\scalebox{1.5}{$,$}}
    \put(50,0){\scalebox{0.2}{\includegraphics{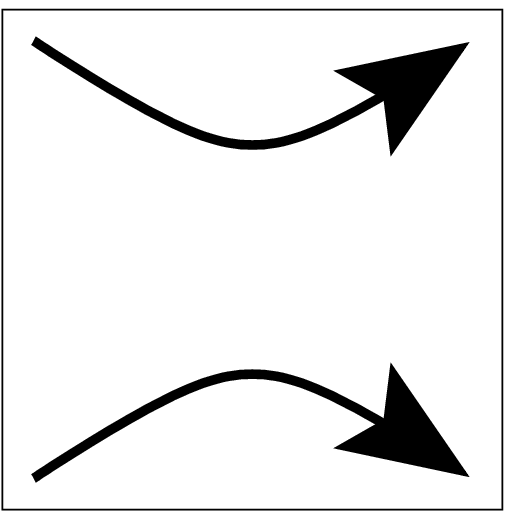}}}
    \put(86,8){\scalebox{1.5}[2.8]{$\left.\right\}$}}
   }
   \put(210,80){
    \put(-43,13){$\Dzero=$}
    \put(-12,8){\scalebox{1.5}[2.8]{$\left\{\right.$}}
    \put( 3,0){\scalebox{0.2}{\includegraphics{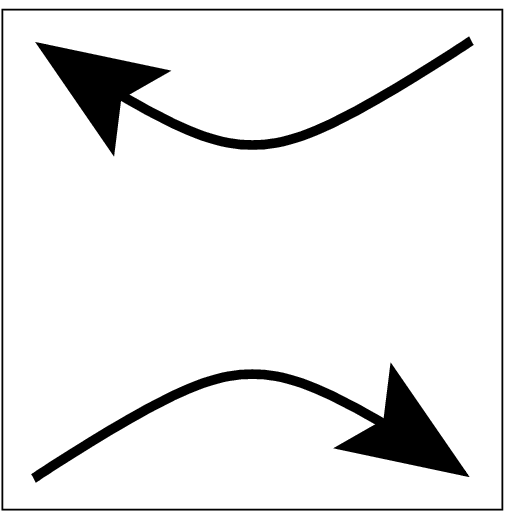}}}
    \put(41,3){\scalebox{1.5}{$,$}}
    \put(53,0){\scalebox{0.2}{\includegraphics{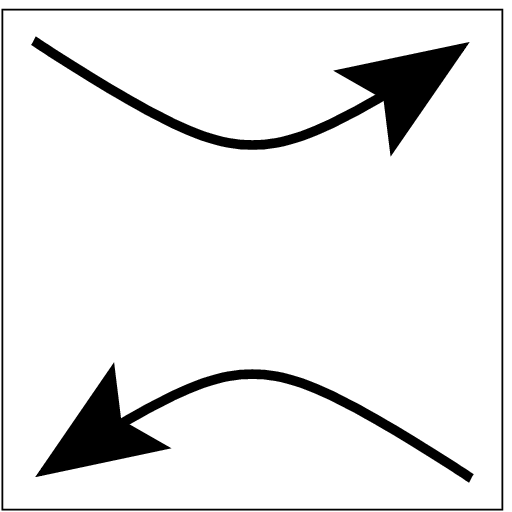}}}
    \put(89,8){\scalebox{1.5}[2.8]{$\left.\right\}$}}
   }
   \put( 50,40){
    \put(-43,13){$\Vinf=$}
    \put(-15,8){\scalebox{1.5}[2.8]{$\left\{\right.$}}
    \put( 0,0){\scalebox{0.2}{\includegraphics{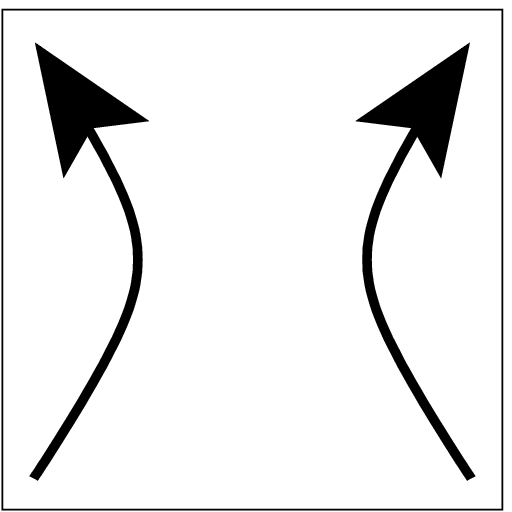}}}
    \put(38,3){\scalebox{1.5}{$,$}}
    \put(50,0){\scalebox{0.2}{\includegraphics{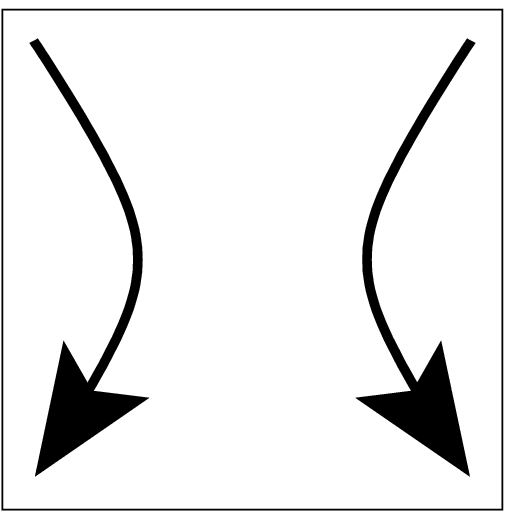}}}
    \put(86,8){\scalebox{1.5}[2.8]{$\left.\right\}$}}
   }
   \put(210,40){
    \put(-43,13){$\Dinf=$}
    \put(-12,8){\scalebox{1.5}[2.8]{$\left\{\right.$}}
    \put( 3,0){\scalebox{0.2}{\includegraphics{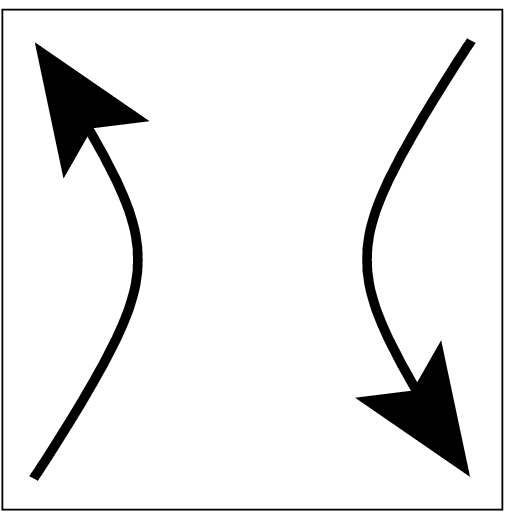}}}
    \put(41,3){\scalebox{1.5}{$,$}}
    \put(53,0){\scalebox{0.2}{\includegraphics{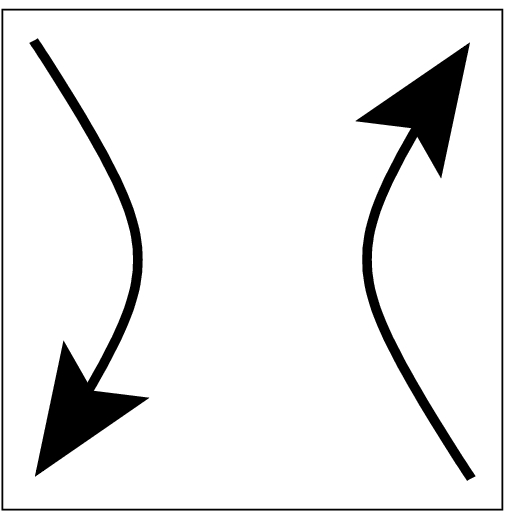}}}
    \put(89,8){\scalebox{1.5}[2.8]{$\left.\right\}$}}
   }
   \put( 50, 0){
    \put(-43,13){$\Vone=$}
    \put(-15,8){\scalebox{1.5}[2.8]{$\left\{\right.$}}
    \put( 0,0){\scalebox{0.2}{\includegraphics{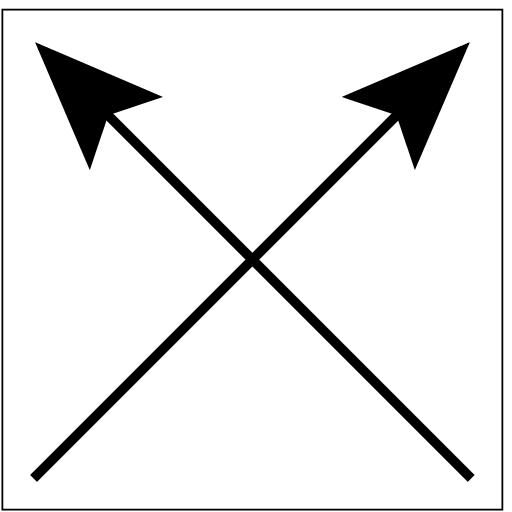}}}
    \put(38,3){\scalebox{1.5}{$,$}}
    \put(50,0){\scalebox{0.2}{\includegraphics{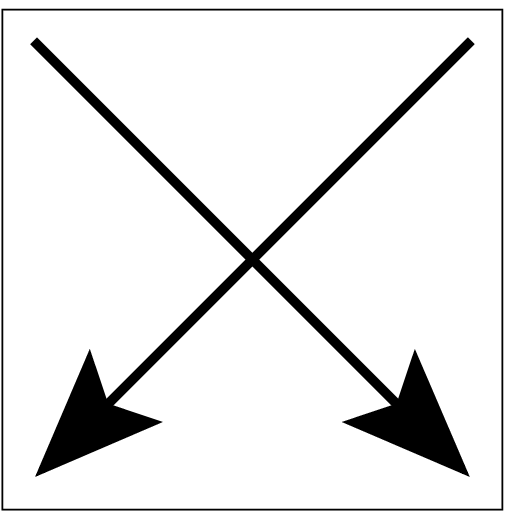}}}
    \put(86,8){\scalebox{1.5}[2.8]{$\left.\right\}$}}
   }
   \put(210, 0){
    \put(-43,13){$\Hone=$}
    \put(-12,8){\scalebox{1.5}[2.8]{$\left\{\right.$}}
    \put( 03,0){\scalebox{0.2}{\includegraphics{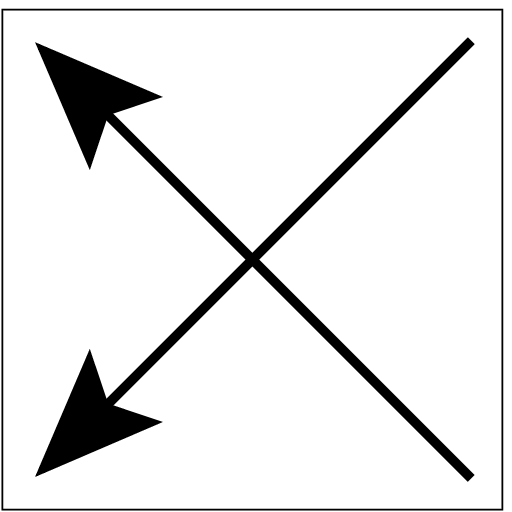}}}
    \put(41,3){\scalebox{1.5}{$,$}}
    \put(53,0){\scalebox{0.2}{\includegraphics{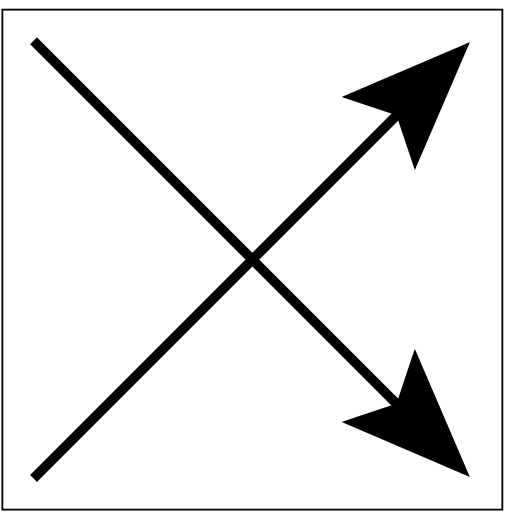}}}
    \put(89,8){\scalebox{1.5}[2.8]{$\left.\right\}$}}
   }
  \end{picture}
 \end{center}
 \caption{$6$ types of oriented tangles}
 \label{Fig:OrientedTangleDiagram}
\end{figure}

When we identify two oriented tangles by homotopy,
all the oriented tangles are divided into $6$ classes
as in Figure \ref{Fig:OrientedTangleDiagram}.
We may call the class including one oriented tangle $\oriented{T}$
the {\em type of the oriented tangle $\oriented{T}$}.
Six types are denoted by
$\Hzero$ (horizontal $0/1$),
$\Dzero$ (diagonal $0/1$),
$\Vinf$ (vertical $1/0$),
$\Dinf$ (diagonal $1/0$),
$\Vone$ (vertical $1/1$),
and $\Hone$ (horizontal $1/1$).
$H_*$ denotes the union of two classes $\Hzero$ and $\Hone$.
Similarly,
we define $V_*$ and $D_*$.

\noindent{\bf Reduced expression of edgepath} \\
We introduce notations
for reduced expressions of edgepaths,
which are used when one edgepath is defined with another edgepath.
%
Let $\Edgepath$ be 
an edgepath
$\angleb{p_k/q_k}$\,--\,%
$\angleb{p_{k-1}/q_{k-1}}$\,--\,%
$\ldots$\,--\,%
$\angleb{p_2/q_2}$\,--\,%
$\angleb{p_1/q_1}$.
Then,
for an integer $a$,
we introduce reduced notations 
of edgepaths:
\begin{eqnarray*}
\lbrack a+\Edgepath\rbrack&=&
\textrm{
$\angleb{a+p_k/q_k}$\,--\,%
$\angleb{a+p_{k-1}/q_{k-1}}$\,--\,%
$\ldots$\,--\,%
$\angleb{a+p_2/q_2}$\,--\,%
$\angleb{a+p_1/q_1}$%
},\\
\lbrack 1/(a+\Edgepath)\rbrack&=&
\textrm{
$\angleb{1/(a+p_k/q_k)}$\,--\,%
$\angleb{1/(a+p_{k-1}/q_{k-1})}$
}\\&&\ \ \ \ 
\textrm{
\,--\,%
$\ldots$\,--\,%
$\angleb{1/(a+p_2/q_2)}$\,--\,%
$\angleb{1/(a+p_1/q_1)}$%
},\\
\lbrack -\Edgepath\rbrack&=&
\textrm{
$\angleb{-p_k/q_k}$\,--\,%
$\angleb{-p_{k-1}/q_{k-1}}$\,--\,%
$\ldots$\,--\,%
$\angleb{-p_2/q_2}$\,--\,%
$\angleb{-p_1/q_1}$%
}.
\end{eqnarray*}

\noindent {\bf Edgepath $\NearlyBasicEdgepath_s(\oriented{T})$} \\
We define $\EdgepathSystem_s(K)$ step by step.
We first define an edgepath $\NearlyBasicEdgepath_s(\oriented{T})$
for an oriented rational tangle $\oriented{T}$ corresponding to $0<\TangleP/\TangleQ<1$.

First,
assume that $\oriented{T}$ corresponds to $\TangleP/\TangleQ=[0,a]=1/a$ with $a\ge 2$.
Let $\oriented{D}$ be the standard diagram of $\oriented{T}$.
If the sign of crossings in $\oriented{D}$ is positive,
we define $\NearlyBasicEdgepath_s(\oriented{T})$ to be
an edgepath
$\angleb{1}$\,--\,$\angleb{1/2}$\,--\,%
$\ldots$\,--\,$\angleb{1/a}$.
%
If the sign of crossings in $\oriented{D}$ is negative,
we define $\NearlyBasicEdgepath_s(\oriented{T})$ to be
an edgepath
$\angleb{0}$\,--\,$\angleb{1/a}$.

  \begin{figure}[hbt]
   \begin{center}
    \begin{picture}(350,70)
     \put(  0,0){
      \put(0,0){\scalebox{0.5}{\includegraphics{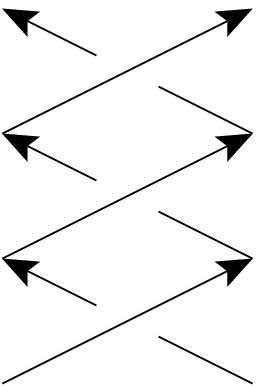}}}
      \put(14,1){$+$}
      \put(14,20){$+$}
      \put(14,39){$+$}
     }
     \put(95,0){
      \put(0,0){\scalebox{0.4}{\includegraphics{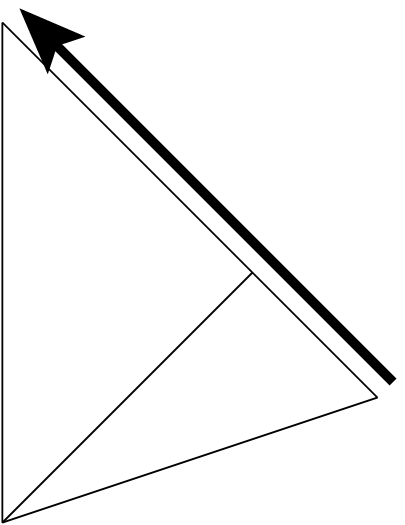}}}
      \put(26,41){$\NearlyBasicEdgepath_s(\oriented{T})$}
      \put(-16,-3){$\angleb{0}$}
      \put(-16,60){$\angleb{1}$}
      \put(48,14){$\angleb{1/a}$}
     }
     \put(190,0){
      \put(0,0){\scalebox{0.5}{\includegraphics{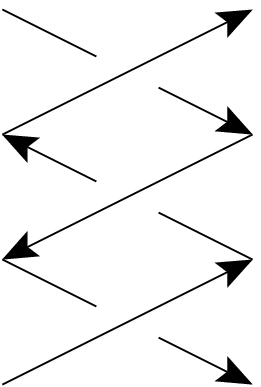}}}
      \put(14,1){$-$}
      \put(14,20){$-$}
      \put(14,39){$-$}
     }
     \put(285,0){
      \put(0,0){\scalebox{0.4}{\includegraphics{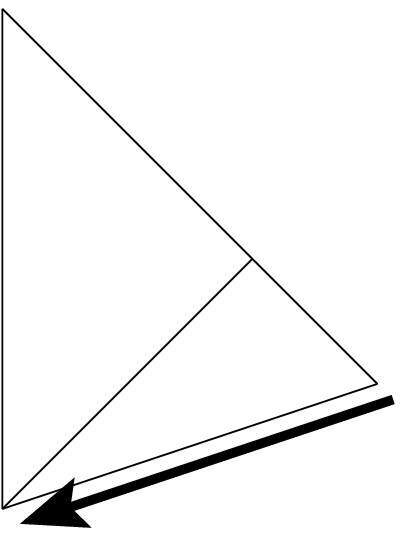}}}
      \put(26,-2){$\NearlyBasicEdgepath_s(\oriented{T})$}
      \put(-16,-3){$\angleb{0}$}
      \put(-16,60){$\angleb{1}$}
      \put(48,17){$\angleb{1/a}$}
     }
    \end{picture}
   \end{center}
   \caption{Diagrams and edgepaths $\NearlyBasicEdgepath_s$ for $1/a$-tangles.}
  \end{figure}

\begin{figure}[hbt]
        \begin{center}
         \begin{picture}(300,80)
          \put(  0,0){
           \put(0,0){\scalebox{0.5}{\includegraphics{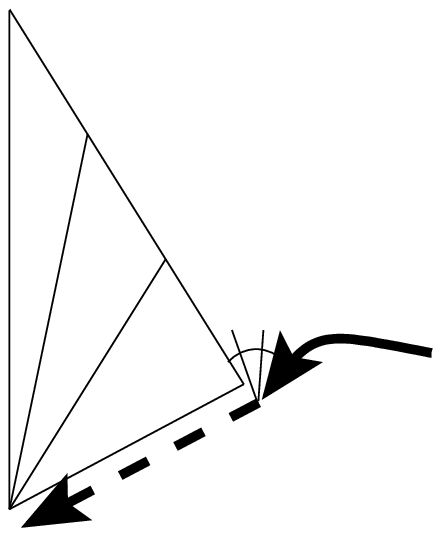}}}
           \put(68,33){odd}
           \put(68,10){$\angleb{1/(a+1)}$}
          }
          \put(100,0){
           \put(0,0){\scalebox{0.5}{\includegraphics{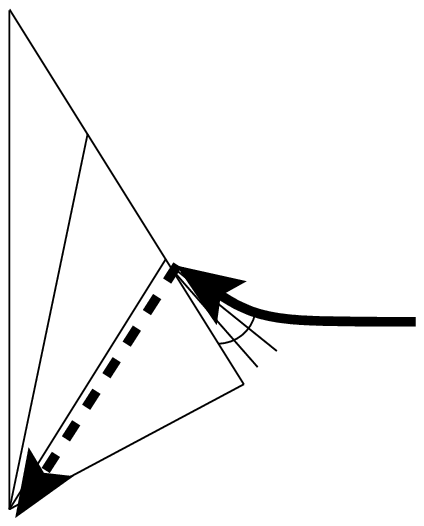}}}
           \put(73,18){odd}
           \put(60,42){$\angleb{1/a}$}
          }
          \put(200,0){
           \put(0,0){\scalebox{0.5}{\includegraphics{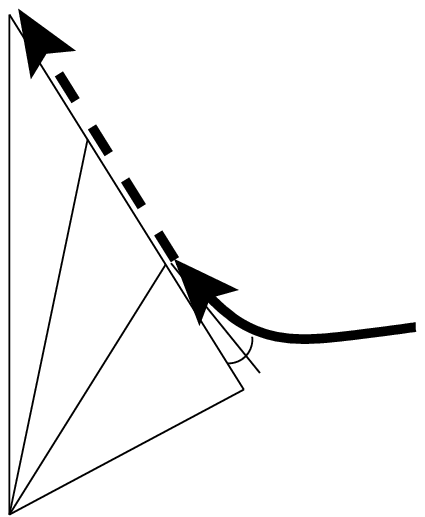}}}
           \put(73,16){even}
           \put(60,42){$\angleb{1/a}$}
          }
         \end{picture}
        \end{center}
        \caption{$1/(a+\NearlyBasicEdgepath_s(\oriented{T_A}))$ and
                 $\NearlyBasicEdgepath_s(\oriented{T_B})$
                 for $\NearlyBasicEdgepath_s(\oriented{T_A})$
                 of type $V_*$, $D_*$ and $H_*$}
        \label{Fig:EdgepathS}
       \end{figure}

Next, 
as an induction hypothesis,
we assume that
$\NearlyBasicEdgepath_s(\oriented{T_A})$ is defined for 
any $\oriented{T_A}$ with the standard continued fraction $[0,a_2,a_3,\ldots,a_k]$ of length $k$ for some $k\ge 2$.
%
For $\oriented{T_A}$ and an integer $a\ge 1$,
let $\oriented{T_B}$ be an oriented tangle
corresponding to $[0,a,a_2,a_3,\ldots,a_k]$ of length $k+1$.
The orientation of $\oriented{T_B}$ is naturally determined by 
the orientation of $\oriented{T_A}$.
$\NearlyBasicEdgepath_s(\oriented{T_B})$ is defined
by using $1/(a+\NearlyBasicEdgepath_s(\oriented{T_A}))$.
\begin{itemize}
\item[(1)]
 Assume that $\oriented{T_A}$ is of type $\Vinf$ or $\Vone$.
$\NearlyBasicEdgepath_s(\oriented{T_A})$ is an edgepath
$\angleb{p_n/q_n}$\,--\,%
$\ldots$\,--\,%
$\angleb{p_1/q_1}$
with $p_n/q_n=1$.
Then,
we define $\NearlyBasicEdgepath_s(\oriented{T_B})$ to be 
an edgepath
$\angleb{0}$\,--\,$\lbrack 1/(a+\NearlyBasicEdgepath_s(\oriented{T_A}))\rbrack$%
$=$%
$\angleb{0}$\,--\,%
$\angleb{1/(a+p_n/q_n)}$\,--\,%
$\ldots$\,--\,%
$\angleb{1/(a+p_1/q_1)}$,
where $1/(a+p_n/q_n)=1/(a+1)$.
\item[(2)]
 Assume that $\oriented{T_A}$ is of type $\Dinf$ or $\Dzero$.
$\NearlyBasicEdgepath_s(\oriented{T_A})$ is an edgepath
$\angleb{p_n/q_n}$\,--\,%
$\ldots$\,--\,%
$\angleb{p_1/q_1}$
with $p_n/q_n=0$.
Then, 
we define $\NearlyBasicEdgepath_s(\oriented{T_B})$ to be 
an edgepath
$\angleb{0}$\,--\,$\lbrack 1/(a+\NearlyBasicEdgepath_s(\oriented{T_A}))\rbrack$%
$=$%
$\angleb{0}$\,--\,%
$\angleb{1/(a+p_n/q_n)}$\,--\,%
$\ldots$\,--\,%
$\angleb{1/(a+p_1/q_1)}$,
where $1/(a+p_n/q_n)=1/a$.
\item[(3)]
 Assume that $\oriented{T_A}$ is of type $\Hzero$ or $\Hone$.
$\NearlyBasicEdgepath_s(\oriented{T_A})$ is an edgepath
$\angleb{p_n/q_n}$\,--\,%
$\ldots$\,--\,%
$\angleb{p_1/q_1}$
with $p_n/q_n=0$.
Then, 
we define $\NearlyBasicEdgepath_s(\oriented{T_B})$ to be 
an edgepath
$\angleb{1}$\,--\,%
$\angleb{1/2}$\,--\,%
$\ldots$\,--\,%
$\angleb{1/(a-1)}$\,--\,%
$\lbrack 1/(a+\NearlyBasicEdgepath_s(\oriented{T_A}))\rbrack$
=
$\angleb{1}$\,--\,%
$\angleb{1/2}$\,--\,%
$\ldots$\,--\,%
$\angleb{1/(a-1)}$\,--\,%
$\angleb{1/(a+p_n/q_n)}$\,--\,%
$\ldots$\,--\,%
$\angleb{1/(a+p_1/q_1)}$,
where $1/(a+p_n/q_n)=1/a$.
\end{itemize}
In Figure \ref{Fig:EdgepathS},
an edgepath $1/(a+\BasicEdgepath_s^\prime(\oriented{T_A}))$ is depicted by a solid curve,
added edges are shown by a dotted line,
and $\BasicEdgepath_s^\prime(\oriented{T_B})$ is
expressed by the combination of a solid curve and a dotted line.
$\NearlyBasicEdgepath_s$ has been defined inductively.
\\
%
%
%

\noindent {\bf Edgepath $\Edgepath_s(\oriented{T})$} \\
Next, 
with the edgepath $\NearlyBasicEdgepath_s$,
we define the edgepath $\Edgepath_s(\oriented{T})$
for an oriented tangle $\oriented{T}$ corresponding to a non-integral $\TangleP/\TangleQ$

First,
we think about an oriented $\TangleR/\TangleS$-tangle
with positive non-integral $\TangleR/\TangleS$.
Assume that 
$z$ is the integer part of $\TangleR/\TangleS$
and $\TangleP/\TangleQ$ is the fractional part,
that is,
$\TangleR/\TangleS=z+\TangleP/\TangleQ$
where
$z\ge 0$ is an integer
and $\TangleP/\TangleQ$ satisfies $0<\TangleP/\TangleQ<1$.
Let $\oriented{T_A}$ be an oriented $\TangleP/\TangleQ$-tangle.
Let $\oriented{T_B}$ be an oriented $\TangleR/\TangleS$-tangle
whose orientation is derived from $\oriented{T_A}$.
%
%
We define $\Edgepath_s(\oriented{T_B})$
as a type III edgepath 
$\angleb{\infty}$\,--\,$\lbrack z+\NearlyBasicEdgepath_s(\oriented{T_A})\rbrack$
if $\oriented{T_A}$ is of type $V_*$ or $D_*$,
and 
as a type II edgepath 
$\angleb{0}$\,--\,$\angleb{1}$\,--\,$\ldots$\,--\,$\angleb{z-1}$%
\,--\,$\lbrack z+\NearlyBasicEdgepath_s(\oriented{T_A})\rbrack$
if $\oriented{T_A}$ is of type $H_*$.

Next, we think about an oriented $\TangleR/\TangleS$-tangle
with negative non-integral $\TangleR/\TangleS$.
Let $\oriented{T_A}$ be an $\TangleR/\TangleS$-tangle with negative $\TangleR/\TangleS$.
Let $\oriented{T_B}$ be an $|\TangleR/\TangleS|$-tangle
obtained from $\oriented{T_A}$
by taking mirror image with respect to the paper plane.
%
We define
$\Edgepath_s(\oriented{T_A})
=\lbrack -\Edgepath_s(\oriented{T_B})\rbrack$.
%

%
%
%
\noindent {\bf Edgepath system $\EdgepathSystem_s(K)$} \\
Last, we define the edgepath system $\EdgepathSystem_s(K)$
for a Montesinos knot $K$.

We give an orientation to a Montesinos knot $K$.
Though there are two choices of orientation,
they will give the same consequence eventually.
As an orientation is given to $K$,
an orientation is given also to each rational tangle $T_i$ of $K$,
so let $\oriented{T_i}$ denote the oriented tangle.

For each oriented tangle $\oriented{T_i}$,
an edgepath $\Edgepath_{s,i}=\Edgepath_s(\oriented{T_i})$ is uniquely determined
as described above. 
%
Let $\EdgepathSystem_s$ denote the edgepath system
obtained by collecting $\Edgepath_{s,i}$ for each $\oriented{T_i}$
of a Montesinos knot $K$.

In \cite{HO},
an edgepath system of a Seifert surface is constructed
as an edgepath system satisfying ``turning number condition'' and ``penultimate vertex condition''.
We can confirm that 
$\EdgepathSystem_s(K)$ satisfies these conditions
and coincides with
the edgepath system of a Seifert surface described in \cite{HO}.
Hence, we have the following.

\begin{lemma}
\label{Lem:SeifertEdgepathSystem}
For a Montesinos knot $K$,
the edgepath system $\EdgepathSystem_s(K)$ is the edgepath system of a Seifert surface.
\end{lemma}


\subsection{Boundary slopes and the numbers of positive/negative crossings}\label{subsection33}

In this part,
we prove the following lemma,
which relates numerical boundary slopes and the numbers of positive/negative crossings.

%
%
\begin{lemma}
 \label{Lem:IdentitiesForKnot}
Let $K$ be a Montesinos knot.
Let 
$\EdgepathSystem_\inc$, 
$\EdgepathSystem_\dec$ and
$\EdgepathSystem_s$
denote
the monotonically increasing edgepath system,
the monotonically decreasing edgepath system and
the edgepath system of a Seifert surface defined above
for $K$.
$D$ is the standard diagram of $K$.
Then, the following identities hold:
\begin{eqnarray*}
 &&\Slope(\EdgepathSystem_\inc)
 (=\Twist(\EdgepathSystem_\inc)
    -\Twist(\EdgepathSystem_s))
    =-2\,\Crossing_{-}(D),
 \\
 &&\Slope(\EdgepathSystem_\dec)
 (=\Twist(\EdgepathSystem_\dec)
     -\Twist(\EdgepathSystem_s))
     =2\,\Crossing_{+}(D)
     .
\end{eqnarray*}
\end{lemma}

\begin{proof}
Recall that
edgepath systems
$\EdgepathSystem_\inc$, $\EdgepathSystem_\dec$ and $\EdgepathSystem_s$
are defined by
$\EdgepathSystem_\inc=(\Edgepath_{\inc,i})$,
$\EdgepathSystem_\dec=(\Edgepath_{\dec,i})$ and
$\EdgepathSystem_s=(\Edgepath_{s,i})$.
The twist of an edgepath system
is the sum of the twists of all edgepaths in the edgepath system.

On the other hand,
the oriented knot diagram $\oriented{D}$ of the standard diagram $D$
is a combination of standard tangle diagrams $\oriented{D_1}$, $\oriented{D_2}$, $\ldots$, $\oriented{D_\NumTangles}$.
Hence, the number of positive/negative crossings of $\oriented{D}$
is the sum of numbers of positive/negative crossings of each $\oriented{D_i}$.

We first prepare the following claim. 


%
%
\begin{claim}
 \label{Claim:0<P/Q<1}
 Assume that $\oriented{T}$ is an oriented $\TangleP/\TangleQ$-tangle with $0<\TangleP/\TangleQ<1$.
 $T$ denotes the unoriented $\TangleP/\TangleQ$-tangle.
 Let $\BasicEdgepath_\dec(\notoriented{T})$
  and $\BasicEdgepath_\inc(\notoriented{T})$
  denote the monotonically decreasing basic edgepath
  and the monotonically increasing basic edgepath
  for the tangle $T$.
 $\NearlyBasicEdgepath_s(\oriented{T})$ is an edgepath 
 defined above for an oriented tangle $\oriented{T}$.
 Let $\oriented{D}$ denote the oriented standard diagram of $\oriented{T}$.
 For
 $\BasicEdgepath_\dec(\notoriented{T})$,
 $\BasicEdgepath_\inc(\notoriented{T})$,
 $\NearlyBasicEdgepath_s(\oriented{T})$ and $\oriented{D}$,
 the following identities hold.
 \begin{eqnarray*}
  &&\Twist(\BasicEdgepath_\inc(\notoriented{T})) - \Twist(\NearlyBasicEdgepath_s(\oriented{T}))
    = -2\,\Crossing_{-}(\oriented{D}) ,\nonumber \\
  &&\Twist(\BasicEdgepath_\dec(\notoriented{T}))-\Twist(\NearlyBasicEdgepath_s(\oriented{T}))
    = 2\,\Crossing_{+}(\oriented{D}) , \nonumber 
 \end{eqnarray*}
  where
   $\Crossing_{+}(\oriented{D})$
   and $\Crossing_{-}(\oriented{D})$
   denote
   the number of positive/negative crossings in $\oriented{D}$
   respectively.
\end{claim}

These are shown by following the definition of each edgepath and edgepath system
and by comparing twists and the numbers of crossings.
\begin{proof}[Proof of Claim \ref{Claim:0<P/Q<1}]
%
%
First, we prove this claim for $1/a$ tangles $(a\ge 2)$.
The standard diagram and three edgepaths are
shown 
for each of two oriented $1/a$-tangles
in Figure \ref{Fig:DiagramsEdgepaths1/aTangles}.

%

\begin{itemize}
 \item
  In the case $\oriented{T}$ has positive crossings,
  edgepaths are
    \begin{eqnarray*}
       && \textrm{
            $\BasicEdgepath_\inc(\notoriented{T})
            =\angleb{1}$\,--\,$\angleb{1/2}$\,--\,$\ldots$\,--\,$\angleb{1/a}$, 
          } 
          \textrm{
            $\BasicEdgepath_\dec(\notoriented{T})
            =\angleb{0}$\,--\,$\angleb{1/a}$, 
          } 
      \\&&
          \textrm{
            $\NearlyBasicEdgepath_s(\oriented{T})
            =\angleb{1}$\,--\,$\angleb{1/2}$\,--\,
            $\ldots$\,--\,$\angleb{1/a}$. 
          }
    \end{eqnarray*}
  The number of crossings are
    \begin{eqnarray*}
      && \Crossing_{+}(\oriented{D})=a, 
      ~~ \Crossing_{-}(\oriented{D})=0.
    \end{eqnarray*}
  Twists are
    \begin{eqnarray*}
      && \Twist(\BasicEdgepath_\inc(\notoriented{T}))=-2(a-1), ~~
         \Twist(\BasicEdgepath_\dec(\notoriented{T}))=2, ~~ \\
      && \Twist(\NearlyBasicEdgepath_s(\oriented{T}))=-2(a-1).
    \end{eqnarray*}
  Hence, we have
    \begin{eqnarray*}
      && \Twist(\BasicEdgepath_\inc(\notoriented{T}))-\Twist(\NearlyBasicEdgepath_s(\oriented{T}))
      =0=-2\,\Crossing_{-}(\oriented{D}), \\
      && \Twist(\BasicEdgepath_\dec(\notoriented{T}))-\Twist(\NearlyBasicEdgepath_s(\oriented{T}))
      =2a=2\,\Crossing_{+}(\oriented{D}).
    \end{eqnarray*}
 \item
 On the other hand,
  in the case $\oriented{T}$ has negative crossings,
  we have: 
    \begin{eqnarray*}
      && \textrm{
            $\BasicEdgepath_\inc(\notoriented{T})
            =\angleb{1}$\,--\,$\angleb{1/2}$\,--\,$\ldots$\,--\,$\angleb{1/a}$,
         } ~
        \textrm{
            $\BasicEdgepath_\dec(\notoriented{T})
            =\angleb{0}$\,--\,$\angleb{1/a}$,
         } ~ \\
         &&
         \textrm{
            $\NearlyBasicEdgepath_s(\oriented{T})
            =\angleb{0}$\,--\,$\angleb{1/a}$. 
         }
    \end{eqnarray*}
  The numbers of crossings are
    \begin{eqnarray*}
      && \Crossing_{+}(\oriented{D})=0, 
      ~~ \Crossing_{-}(\oriented{D})=a. 
    \end{eqnarray*}
  Twists are
    \begin{eqnarray*}
      && \Twist(\BasicEdgepath_\inc(\notoriented{T}))=-2(a-1),  
      ~~ \Twist(\BasicEdgepath_\dec(\notoriented{T}))=2, 
      ~~ \Twist(\NearlyBasicEdgepath_s(\oriented{T}))=2.
    \end{eqnarray*}
  Hence, we have
    \begin{eqnarray*}
      && \Twist(\BasicEdgepath_\inc(\notoriented{T}))-\Twist(\NearlyBasicEdgepath_s(\oriented{T}))
      = -2a=-2\,\Crossing_{-}(\oriented{D}), \\
      && \Twist(\BasicEdgepath_\dec(\notoriented{T}))-\Twist(\NearlyBasicEdgepath_s(\oriented{T}))
      = 0=2\,\Crossing_{+}(\oriented{D}).
    \end{eqnarray*}
\end{itemize}
Eventually, we have proved the claim
for $T$ with the standard continued fraction $[0,a]$ of length $2$.
  \begin{figure}[hbt]
   \begin{center}
    \begin{picture}(350,70)
     \put(  0,0){
      \put(0,0){\scalebox{0.5}{\includegraphics{positive-one-over-odd.eps}}}
      \put(14,1){$+$}
      \put(14,20){$+$}
      \put(14,39){$+$}
     }
     \put(95,0){
      \put(0,0){\scalebox{0.4}{\includegraphics{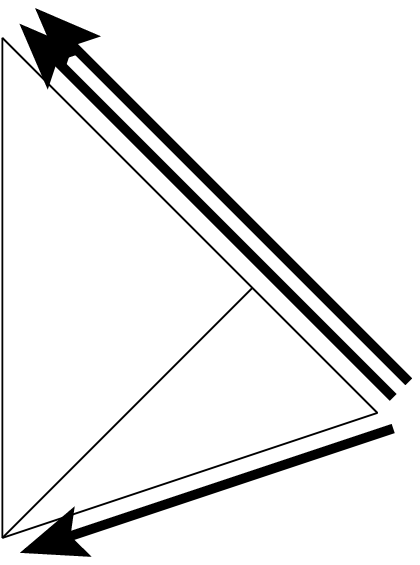}}}
      \put(-15,-3){$\angleb{0}$}
      \put(-16,60){$\angleb{1}$}
      \put(48,16){$\angleb{1/a}$}
      \put(-5,25){\vector(2,1){30}}
      \put(-39,23){$\NearlyBasicEdgepath_s(\oriented{T})$}
      \put(26,0){$\BasicEdgepath_\dec(\notoriented{T})$}
      \put(26,46){$\BasicEdgepath_\inc(\notoriented{T})$}
     }
     \put(190,0){
      \put(0,0){\scalebox{0.5}{\includegraphics{negative-one-over-odd.eps}}}
      \put(14,1){$-$}
      \put(14,20){$-$}
      \put(14,39){$-$}
     }
     \put(285,0){
      \put(0,0){\scalebox{0.4}{\includegraphics{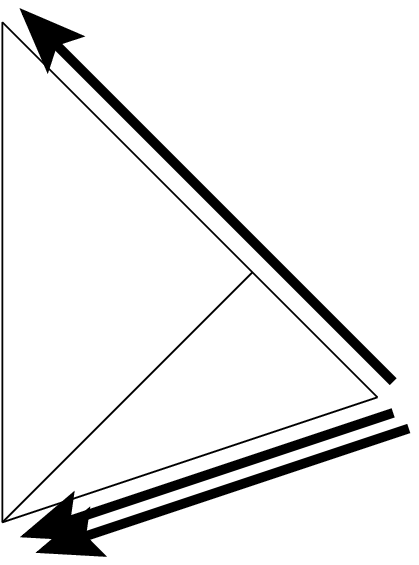}}}
      \put(-15,-3){$\angleb{0}$}
      \put(-16,60){$\angleb{1}$}
      \put(48,16){$\angleb{1/a}$}
      \put(-8,26){\vector(2,-1){30}}
      \put(-42,23){$\NearlyBasicEdgepath_s(\oriented{T})$}
      \put(26,0){$\BasicEdgepath_\dec(\notoriented{T})$}
      \put(26,44){$\BasicEdgepath_\inc(\notoriented{T})$}
     }
    \end{picture}
   \end{center}
   \caption{Diagrams and edgepaths for $1/a$-tangles.}
   \label{Fig:DiagramsEdgepaths1/aTangles}
  \end{figure}

%
%

Next,
as a hypothesis of induction,
assume that the claim is true
for $\oriented{T_A}$ with the standard continued fraction $[0,a_2,\ldots,a_{k}]$ of length $k$ for fixed $k\ge 2$.
Hence, 
we have
 \begin{eqnarray*}
  \Twist(\BasicEdgepath_\inc(\notoriented{T_A}))-\Twist(\NearlyBasicEdgepath_s(\oriented{T_A}))
    &=& -2\,\Crossing_{-}(\oriented{D_A}) ,\\
 \Twist(\BasicEdgepath_\dec(\notoriented{T_A}))-\Twist(\NearlyBasicEdgepath_s(\oriented{T_A}))
    &=& 2\,\Crossing_{+}(\oriented{D_A}) .
 \end{eqnarray*}
Assume that $0<\TangleP/\TangleQ<1$ and $\TangleP/\TangleQ=[0,a_2,\ldots,a_k]$.
Let $T_A$ be a rational tangle corresponding to $\TangleP/\TangleQ$. 
Let $\oriented{T_A}$ be a rational tangle $T_A$ with either of orientation.
Assume that $a\ge 1$.
Let $\oriented{T_B}$ be a rational tangle with orientation
obtained by
first taking mirror image with respect to a line from upper left to lower right
and then adding a $1/a$-tangle below.
Note that orientation of $\oriented{T_B}$ is naturally given 
by the orientation of $\oriented{T_A}$.
$T_B$ is a rational tangle corresponding to the fraction 
$1/(a+\TangleP/\TangleQ)=[0,a,a_1,a_2,\ldots,a_k]$,
which is denoted by $\TangleR/\TangleS$.

       %
       %
       Assume that $\oriented{T_A}$ is of type $\Vinf$ or $\Vone$.
       The oriented standard diagram $\oriented{D_A}$ of $\oriented{T_A}$
       and
       the oriented standard diagram $\oriented{D_B}$ of $\oriented{T_B}$
       are shown in the top two figures 
       in Figure \ref{Fig:DiagramsEdgepaths0<p/q<1ForV*}.
       In the bottom two figures,
       $\NearlyBasicEdgepath_s$,
       $\BasicEdgepath_\dec$
       and
       $\BasicEdgepath_\inc$ are depicted
       for $\oriented{T_A}$ and $\oriented{T_B}$.
       In the bottom right figure,
       a dotted segment represents added edges,
       a solid curve is an edgepath obtained from an edgepath for $\oriented{T_A}$
       and 
       each of
       $\BasicEdgepath_\inc(\oriented{T_B})$,
       $\BasicEdgepath_\dec(\oriented{T_B})$ and
       $\NearlyBasicEdgepath_s(\oriented{T_B})$
       corresponds to 
       the combination of a dotted edge and a solid curve.
%
 %
%
%
       \begin{figure}[hbt]
        \begin{center}
         \begin{picture}(250,200)
          \put( 10,130){
           \put(0,0){\scalebox{0.3}{\includegraphics{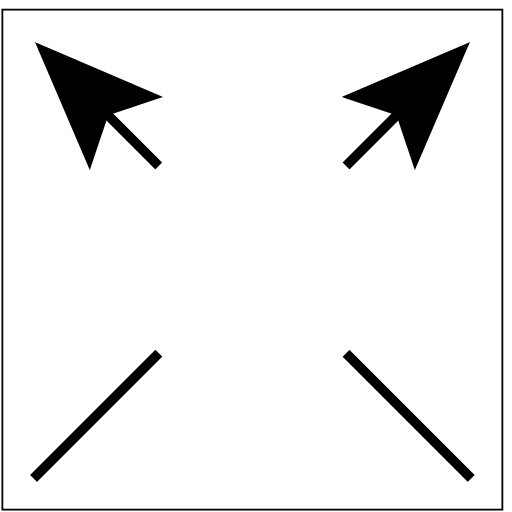}}}
           \put(16,18){$\oriented{D_A}$}
          }
          \put(155,110){
           \put(0,0){\scalebox{0.3}{\includegraphics{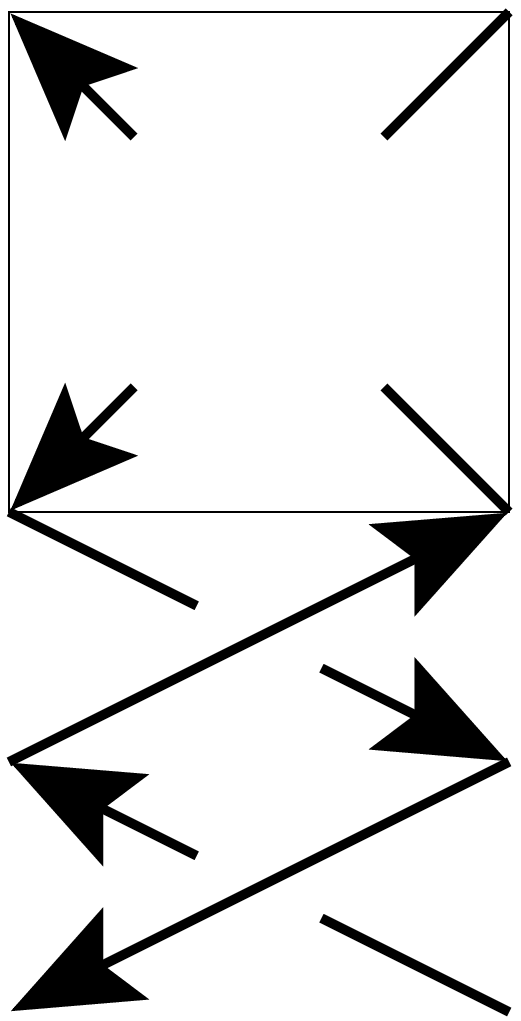}}}
           \put(16,60){\rotatebox{90}{\reflectbox{$\oriented{D_A}$}}}
           \put(5,33){$-$}
           \put(32,10){$-$}
           \put(48,43){\rotatebox{-90}{\scalebox{1.0}[1.0]{$\overbrace{\hspace{15mm}}$}}}
           \put(58,20){$a$}
          }
          \put(0,0){
           \put(0,0){\VDiagram}
           \put(47,23){$\angleb{\TangleP/\TangleQ}$}
          }
          \put(100,0){
           \put(0,0){\scalebox{0.5}{\includegraphics{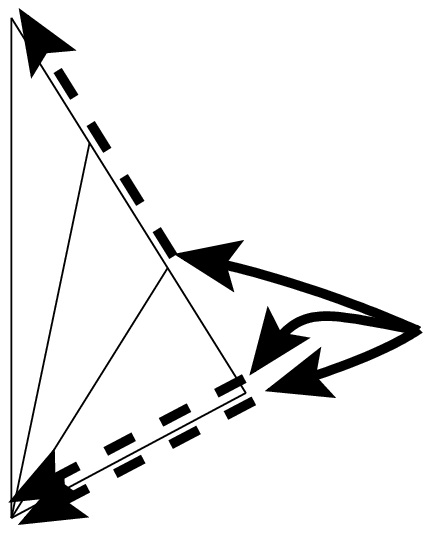}}}
           \put(28,81){$\angleb{1}$}
           \put(28,-8){$\angleb{0}$}
           \put(78,3){$1/(a+\BasicEdgepath_\inc(\notoriented{T_A}))$}
           \put(88,12){\vector(0,1){11}}
           \put(79,54){$1/(a+\BasicEdgepath_\dec(\notoriented{T_A}))$}
           \put(85,50){\vector(0,-1){13}}
           \put(49,81){$1/(a+\NearlyBasicEdgepath_s(\oriented{T_A}))$}
           \put(77,77){\vector(0,-1){45}}
          \put(28,39){\vector(1,0){30}}
          \put(0,36){$\angleb{1/a}$}
           \put(72,0){\vector(0,1){17}}
           \put(60,-10){$\angleb{1/(a+1)}$}
           \put(97,28){$\angleb{\TangleR/\TangleS}$}
          }
         \end{picture}
        \end{center}
        \caption{$\oriented{D_A}$, $\oriented{D_B}$, 
                 edgepaths for $\oriented{T_A}$, 
                 edgepaths for $\oriented{T_B}$
                 in the case $\oriented{T_A}$ is of $V_*$ }
        \label{Fig:DiagramsEdgepaths0<p/q<1ForV*}
       \end{figure}

       Edgepaths are
         \begin{eqnarray*}
           && \textrm{
                 $\BasicEdgepath_\inc(\notoriented{T_B})
                 =\angleb{1}$\,--\,$\angleb{1/2}$\,--\,$\ldots$%
                 \,--\,$\angleb{1/(a-1)}$%
                 \,--\,$\lbrack 1/(a+\BasicEdgepath_\dec(\notoriented{T_A}))\rbrack$ } \\
          && \hspace{14mm} \textrm{
                 $=\angleb{1}$\,--\,$\angleb{1/2}$\,--\,$\ldots$%
                 \,--\,$\angleb{1/a}$\,--\,$\ldots$
                 \,--\,$\angleb{\TangleR/\TangleS}$,
              }
           \\
           && \textrm{
                 $\BasicEdgepath_\dec(\notoriented{T_B})
                 =\angleb{0}$%
                 \,--\,$\lbrack 1/(a+\BasicEdgepath_\inc(\notoriented{T_A}))\rbrack$%
                 $=\angleb{0}$\,--\,$\angleb{1/(a+1)}$%
                 \,--\,$\ldots$\,--\,$\angleb{\TangleR/\TangleS}$,
              } 
           \\
           && \textrm{
                 $\NearlyBasicEdgepath_s(\oriented{T_B})%
                 =\angleb{0}$\,--\,$\lbrack 1/(a+\NearlyBasicEdgepath_s(\oriented{T_A}))\rbrack$%
                 $=\angleb{0}$\,--\,%
                 $\angleb{1/(a+1)}$\,--\,$\ldots$\,--\,$\angleb{\TangleR/\TangleS}$.
              }
         \end{eqnarray*}
       %
       Signs of crossings invert on taking mirror.
       From the orientation of the tangle $\oriented{T_A}$,
       the sign of crossings of $1/a$ tangle is negative.
       Hence, the numbers of crossings are
         \begin{eqnarray*}
           && \Crossing_{+}(\oriented{D_B})=\Crossing_{-}(\oriented{D_A}), \\
           && \Crossing_{-}(\oriented{D_B})=\Crossing_{+}(\oriented{D_A})+a. 
         \end{eqnarray*}
       On making $[1/(a+\Edgepath)]$ from $\Edgepath$,
       the original edgepath $\Edgepath$ is turned roughly upside down.
       Then, the sign of the effect of an edgepath to the twist inverts.
       Twists are
         \begin{eqnarray*}
           && \Twist(\BasicEdgepath_\inc(\notoriented{T_B}))
              =-2(a-1)-\Twist(\BasicEdgepath_\dec(\notoriented{T_A})), \\
           && \Twist(\BasicEdgepath_\dec(\notoriented{T_B}))
              =2-\Twist(\BasicEdgepath_\inc(\notoriented{T_A})), \\
           && \Twist(\NearlyBasicEdgepath_s(\oriented{T_B}))
              =2-\Twist(\NearlyBasicEdgepath_s(\oriented{T_A})). 
         \end{eqnarray*}
       %
       Hence, we have
         \begin{eqnarray*}
           \Twist(\BasicEdgepath_\inc(\notoriented{T_B}))
             -\Twist(\NearlyBasicEdgepath_s(\oriented{T_B}))
           &=&
             -2a
             -\Twist(\BasicEdgepath_\dec(\notoriented{T_A}))
             +\Twist(\NearlyBasicEdgepath_s(\oriented{T_A}))
           \\
           &=&-2a-2\,\Crossing_{+}(\oriented{D_A})
           =-2\,\Crossing_{-}(\oriented{D_B}), \\
           \Twist(\BasicEdgepath_\dec(\notoriented{T_B}))
             -\Twist(\NearlyBasicEdgepath_s(\oriented{T_B}))
           &=&\Twist(\NearlyBasicEdgepath_s(\oriented{T_A}))
             -\Twist(\BasicEdgepath_\inc(\notoriented{T_A}))
           \\
           &=&2\,\Crossing_{-}(\oriented{D_A})
           =2\,\Crossing_{+}(\oriented{D_B}). 
         \end{eqnarray*}

The inductive conclusion has been proved for $T_A$ of type-$V_*$.
Similarly, the conclusion can be shown if $T_A$ is of type-$D_*$ or type-$H_*$.
Thus, the claim holds for $\oriented{T}$
if the corresponding fraction has the standard continued expansion of length $k+1$.
By induction,
these are enough to show the lemma.
\end{proof}

Now we have the following;

\begin{claim}
 \label{Claim:generalP/Q}
   Let $\oriented{T}$ be the oriented $\TangleP/\TangleQ$-tangle with a non-integral $\TangleP/\TangleQ$.
   $T$ denotes the unoriented $\TangleP/\TangleQ$-tangle.
   $\oriented{D}$ denotes the oriented standard diagram of $\oriented{T}$.

   Then, identities
   \[\Twist(\Edgepath_\inc(\notoriented{T}))
     -\Twist(\Edgepath_s(\oriented{T}))
     =-2\,\Crossing_{-}(\oriented{D})
   \]
  and
   \[\Twist(\Edgepath_\dec(\notoriented{T}))
     -\Twist(\Edgepath_s(\oriented{T}))
     =2\,\Crossing_{+}(\oriented{D})
   \]
  hold.
\end{claim}

\begin{proof}
With Claim \ref{Claim:0<P/Q<1},
we can show Claim \ref{Claim:generalP/Q} similarly to the proof of Claim \ref{Claim:0<P/Q<1}.
First,
in order to show this claim for a non-integral $\TangleP/\TangleQ>0$,
we check the effect on adding integral tangles.
Next,
in order to show this claim for a non-integral $\TangleP/\TangleQ<0$,
we check the effect on reversing the sign of $\TangleP/\TangleQ$.
The standard diagram and edgepaths $\Edgepath_\inc$, $\Edgepath_\dec$ and $\Edgepath_s$ for $\TangleP/\TangleQ<0$
is obtained by taking reflections
of the standard diagram and edgepaths for $|\TangleP/\TangleQ|$.
\end{proof}

By taking sum of identities in Claim \ref{Claim:generalP/Q}, 
Lemma \ref{Lem:IdentitiesForKnot} is proved immediately. 
%
%
\end{proof}

\subsection{Main theorem for an arbitrary standard diagram}\label{subsection34}


\begin{proof}[Proof of the main theorem]
By combining 
Lemma \ref{Lem:TwistInequality} and 
Lemma \ref{Lem:IdentitiesForKnot},
we easily have
the inequality in the statement of the main theorem,
for any standard diagram of a Montesinos knot,
which may not attain the crossing number.
\end{proof}

\subsection{A standard diagram attaining the minimal crossing number}\label{subsection35}

%

In this final part, concerning Remark \ref{Rem:StandardMinimal} 
we confirm that 
there is a standard diagram attaining the minimal crossing number.


A standard diagram is defined not for an isotopy class of a Montesinos knot 
but for an expression of a Montesinos knot.
A Montesinos knot 
$M(\TangleP_1/\TangleQ_1$, $\ldots$, $\TangleP_{i-1}/\TangleQ_{i-1}$, $\TangleP_i/\TangleQ_i$, $\TangleP_{i+1}/\TangleQ_{i+1}$, $\TangleP_{i+2}/\TangleQ_{i+2}$, $\ldots$, $\TangleP_\NumTangles/\TangleQ_\NumTangles)$
can be isotoped to
another Montesinos knot
$M(\TangleP_1/\TangleQ_1$, $\ldots$, $\TangleP_{i-1}/\TangleQ_{i-1}$, $(\TangleP_i/\TangleQ_i)\pm 1$, $(\TangleP_{i+1}/\TangleQ_{i+1})\mp 1$, $\TangleP_{i+2}/\TangleQ_{i+2}$, $\ldots$, $\TangleP_\NumTangles/\TangleQ_\NumTangles)$
where the isotopy is just a rotation of one tangle.
With repeated usage of this isotopy,
a Montesinos knot has infinitely many different expressions as a Montesinos knot. 
%
%
By repeatedly using the above isotopy,
every Montesinos knot can be isotoped into another restricted form,
which corresponds to a Montesinos knot with a restricted expression.
A \textit{restricted expression} of a Montesinos knot
is a tuple of fractions satisfying one of
\begin{itemize}
 \item[(i)] all $\TangleP_i/\TangleQ_i$ is positive,
 \item[(ii)] all $\TangleP_i/\TangleQ_i$ is negative,
 \item[(iii)] $(\TangleP_1/\TangleQ_1,\TangleP_2/\TangleQ_2,\ldots,\TangleP_\NumTangles/\TangleQ_\NumTangles)$
            includes both positive and negative fractions,
            and all $\TangleP_i/\TangleQ_i$ has absolute value less than $1$.
\end{itemize}

The standard diagram of a Montesinos knot with a restricted expression
actually gives a ``reduced Montesinos diagram'' defined in \cite{LT}. 
%
By the following theorem determining the crossing numbers of Montesinos knots,
we see that this standard diagram attains the minimal crossing number.
\begin{theorem*}[{\cite[Theorem 10]{LT}}]
If a link L admits an n-crossing, reduced Montesinos diagram, 
then L cannot be projected with fewer than n crossings.
\end{theorem*}
In short,
any Montesinos knot can be expressed by a restricted expression,
and the expression corresponds to the standard diagram which is a minimal diagram. 
%


\section*{Acknowledgement}

The authors would like to thank Koya Shimokawa for useful conversations.

\end{document}